\documentclass[12pt,reqno]{article}

\usepackage[usenames]{color}
\usepackage{amssymb}
\usepackage{graphicx}
\usepackage{amscd}

\usepackage{amsthm}
\newtheorem{theorem}{Theorem}

\newtheorem{conjecture}[theorem]{Conjecture}

\theoremstyle{definition}

\newtheorem{example}[theorem]{Example}

\usepackage[colorlinks=true,
linkcolor=webgreen, filecolor=webbrown,
citecolor=webgreen]{hyperref}

\definecolor{webgreen}{rgb}{0,.5,0}
\definecolor{webbrown}{rgb}{.6,0,0}

\usepackage{color}

\usepackage{float}

\usepackage{graphics,amsmath,amssymb}
\usepackage{amsfonts}
\usepackage{latexsym}
\usepackage{epsf}

\newcommand{\seqnum}[1]{\href{http://oeis.org/#1}{\underline{#1}}}

\begin{document}

\begin{center}
\vskip 1cm{\LARGE\bf  On sequences with $\{-1,0,1\}$ Hankel transforms}  \vskip 1cm
\large
Paul Barry\\
School of Science\\
Waterford Institute of Technology\\
Ireland\\

\end{center}
\vskip .2 in

\begin{abstract} We study Hankel transforms of sequences, where the transform elements are members of the set $\{-1,0,1\}$. We relate these Hankel transforms to special continued fraction expansions. In particular, we posit a conjecture relating the distribution of non-zero terms in the Hankel transform to the distribution of powers of the variable in the defining continued fractions.
\end{abstract}

\section{Introduction} Given a sequence $a_n$, we denote by $h_n$ the general term of the sequence
with $h_n=|a_{i+j}|_{0 \le i,j \le n}$. The sequence $h_n$ is called the Hankel transform of $a_n$ \cite{Layman}.
For a given sequence $a_n$, it can be shown that the sequence
$$\sum_{k=0}^n \binom{n}{k} r^{n-k} a_k$$ will also have the same Hankel transform. Similarly, if
the sequence $a_n$ has generating function $f(x)$, where $f(0)\ne 0$,  then the sequence with generating function
$\frac{f(x)}{1-rxf(x)}$ will also have same Hankel transform. In both cases, $r \in \mathbb{R}$ is arbitrary.
Thus many sequences may have the same Hankel transform, and therefore the problem of inverting the Hankel transform is not an elementary one. In this note, we show, subject to a deep conjecture, that in one instance, the relationship between a sequence $h_n$ and its pre-image is more easily determined.

Although in the sequel we will exhibit Hankel transforms with ostensibly more than one pre-image, we introduce the notion of element multiplicity to confer a degree of uniqueness in the transform.

In the sequel, we will work with integer sequences. We will refer to some known sequences by their ``$Annnnnn$'' number in the On-Line Encyclopedia of Integer Sequences \cite{SL1}.

\section{A conjecture concerning special continued fractions}

\begin{conjecture} Consider a sequence of natural numbers  $p_n$, where $p_0=1$ and the sequence $b_n$ defined by
 \begin{equation}\label{b_seq} b_n=\sum_{k=0}^{\lfloor \frac{n}{2} \rfloor} p_{n-2k} - \frac{1+(-1)^n}{2} \end{equation} is non-negative and non-decreasing, with $b_n \in \mathbb{N}_0$,  $b_0=0,b_1 \ne 0$.   Then the sequence $a_n$ with generating function expressed as the continued fraction
$$\cfrac{1}{1\pm
\cfrac{x^{p_0}}{1\pm
\cfrac{x^{p_1}}{1\pm
\cfrac{x^{p_2}}{1\pm \cdots}}}}$$
has a Hankel transform consisting solely of the numbers $-1$, $0$, and $1$. Moreover, the non-zero terms occur at locations indexed by the sequence $b_n$.

Conversely, given a sequence of numbers $b_n \in N_0$,  $b_0=0,b_1 \ne 0, b_2, b_3,\ldots$ with $b_i \le b_{i+1}$, then the sequence whose generating function is given by the above continued fraction where the sequence $p_n$ is given by
\begin{displaymath}
\left(\begin{array}{c} p_0 \\ p_1 \\p_2\\ p_3 \\p_4 \\p_5 \\ \vdots \end{array}\right)=
\left(\begin{array}{ccccccc}1 & 0 & 0 & 0 & 0 & 0 & \ldots \\
                            0 & 1 & 0 & 0 & 0 & 0 & \ldots \\
                            0 & 0 & 1 & 0& 0 & 0 & \ldots \\
                            0 & -1 & 0 & 1 & 0 & 0 & \ldots \\
                            0 & 0 & -1&  0 & 1 & 0 & \ldots \\
                            0 & 0  & 0 & -1 & 0 & 1 &\ldots\\ \vdots
& \vdots &
\vdots & \vdots & \vdots & \vdots &
\ddots\end{array}\right)\left(\begin{array}{c} 1 \\ b_1 \\b_2\\ b_3 \\b_4 \\b_5 \\ \vdots \end{array}\right)\end{displaymath}
has a $\{-1,0,1\}$ Hankel transform whose non-zero elements are indexed by the sequence $b_n$.
\end{conjecture}

Equation (\ref{b_seq}) can be visualized in the following way.
\begin{displaymath}
\left(\begin{array}{c} b_0 \\ b_1 \\b_2\\ b_3 \\b_4 \\b_5 \\ \vdots \end{array}\right)=
\left(\begin{array}{ccccccc}0 & 0 & 0 & 0 & 0 & 0 & \ldots \\
                            0 & 1 & 0 & 0 & 0 & 0 & \ldots \\
                            0 & 0 & 1 & 0& 0 & 0 & \ldots \\
                            0 & 1 & 0 & 1 & 0 & 0 & \ldots \\
                            0 & 0 & 1&  0 & 1 & 0 & \ldots \\
                            0 & 1  & 0 & 1 & 0 & 1 &\ldots\\ \vdots
& \vdots &
\vdots & \vdots & \vdots & \vdots &
\ddots\end{array}\right)\left(\begin{array}{c} p_0 \\ p_1 \\p_2\\ p_3 \\p_4 \\p_5 \\ \vdots \end{array}\right),\end{displaymath} where column $k$ of the above matrix (except for column $0$) has generating function $\frac{x^k}{1-x^2}$. The sequence of numbers given by
$\sum_{k=0}^{\lfloor \frac{n}{2} \rfloor} p_{n-2k}$ corresponds to the degrees of the denominator polynomials in the partial continued fractions.

In the case that $b_i=b_{i+1}=\cdots =b_{i+r}$  we say that the corresponding non-zero term occurs with \emph{multiplicity} $r+1$.

We shall call the sequence of powers $\{p_n\}$ in the continued fraction the \emph{CF power sequence} and we shall call the indexing sequence $b_n$ the \emph{Hankel pattern sequence}. It specifies the pattern of occurrence of the non-zero terms in the Hankel transform. If we let $P(x)$ denote the generating function of the sequence $p_n$, and $B(x)$ denote the generating function of $b_n$, then we have the relations
\begin{equation} \label{rel_1} B(x)=\frac{1}{1-x^2}P(x)-\frac{1}{1-x^2},\end{equation} and
\begin{equation} \label{rel_2} P(x)=(1-x^2)B(x)-1.\end{equation}

Note that in the continued fraction above, the $\pm$ indicates that an arbitrary choice of ``$+$'' or ``$-$'' is possible at each stage. Note also that the conjecture is silent on the matter of the distribution of the minus signs in the transform. 

\begin{example} The Jacobsthal numbers $J_n=\frac{2^n}{3}-\frac{(-1)^n}{3}$ \seqnum{A001045} give the sequence of numbers
$$0,1,1,3,5,11,21,43,\ldots.$$ We seek to find a sequence $a_n$ whose Hankel transform is composed of the numbers $-1$, $0$ and $1$, where the distribution of the non-zero terms is governed (or indexed) by the Jacobsthal numbers. Thus we want a Hankel transform of the form
$$\alpha,\beta,0,\gamma,0,\delta,0,0,0,0,0,\epsilon,0,0,0,0,0,0,0,0,0,\zeta,0,0,\ldots$$ where $\alpha, \beta,\ldots \in \{-1,1\}$.
As it is not immediately obvious how to deal with the duplicated $1$ ($J_1=J_2=1$), we shall ignore this for the moment and assume that the sequence is $$0,1,3,5,11,21,43,\ldots.$$
We look for the sequence of powers of $x$ to be used in the defining continued fraction.
\begin{displaymath}\left(\begin{array}{c} p_0 \\ p_1 \\p_2\\ p_3 \\p_4 \\p_5 \\ \vdots \end{array}\right)=
\left(\begin{array}{ccccccc}1 & 0 & 0 & 0 & 0 & 0 & \ldots \\
                            0 & 1 & 0 & 0 & 0 & 0 & \ldots \\
                            0 & 0 & 1 & 0& 0 & 0 & \ldots \\
                            0 & -1 & 0 & 1 & 0 & 0 & \ldots \\
                            0 & 0 & -1&  0 & 1 & 0 & \ldots \\
                            0 & 0  & 0 & -1 & 0 & 1 &\ldots\\ \vdots
& \vdots &
\vdots & \vdots & \vdots & \vdots &
\ddots\end{array}\right)\left(\begin{array}{c} 1 \\ 1 \\3\\ 5 \\11 \\21 \\ \vdots \end{array}\right)=
\left(\begin{array}{c} 1 \\ 1 \\3\\ 4 \\8 \\16 \\ \vdots \end{array}\right)\end{displaymath}

This leads us to consider the sequence $a_n$ with generating function
$$\cfrac{1}{1-
\cfrac{x}{1-
\cfrac{x}{1-
\cfrac{x^3}{1-
\cfrac{x^4}{1-
\cfrac{x^8}{1-
\cfrac{x^{16}}{1-\cdots}}}}}}}.$$
Thus $a_n$ begins
$$1, 1, 2, 4, 8, 17, 36, 76, 161, 342, 726, 1541, 3272, 6948, 14753,  \ldots$$ with the desired Hankel transform
$$1, 1, 0, -1, 0, 1, 0, 0, 0, 0, 0, -1, 0, 0, 0, 0, 0, 0, 0, 0, 0, 1, 0,\ldots.$$
Now let us see what happens when we retain both $J_1$ and $J_2$.
We look for the sequence of powers of $x$ to be used in the defining continued fraction.
\begin{displaymath}\left(\begin{array}{c} p_0 \\ p_1 \\p_2\\ p_3 \\p_4 \\p_5 \\ \vdots \end{array}\right)=
\left(\begin{array}{ccccccc}1 & 0 & 0 & 0 & 0 & 0 & \ldots \\
                            0 & 1 & 0 & 0 & 0 & 0 & \ldots \\
                            0 & 0 & 1 & 0& 0 & 0 & \ldots \\
                            0 & -1 & 0 & 1 & 0 & 0 & \ldots \\
                            0 & 0 & -1&  0 & 1 & 0 & \ldots \\
                            0 & 0  & 0 & -1 & 0 & 1 &\ldots\\ \vdots
& \vdots &
\vdots & \vdots & \vdots & \vdots &
\ddots\end{array}\right)\left(\begin{array}{c} 1 \\ 1 \\1\\ 3 \\5 \\11 \\ \vdots \end{array}\right)=
\left(\begin{array}{c} 1 \\ 1 \\1\\ 2 \\4 \\8 \\ \vdots \end{array}\right)\end{displaymath}

This leads us to consider the new sequence $a_n^*$ with generating function
$$\cfrac{1}{1-
\cfrac{x}{1-
\cfrac{x}{1-
\cfrac{x}{1-
\cfrac{x^2}{1-
\cfrac{x^4}{1-
\cfrac{x^8}{1-\cdots}}}}}}}.$$ We obtain the sequence that starts
$$1, 1, 2, 5, 13, 35, 95, 259, 707, 1932, 5281, 14438, 39475, 107933, 295115, 806922, 2206342, \ldots,$$ and which has Hankel transform
$$1, 1, 0, -1, 0, 1, 0, 0, 0, 0, 0, -1, 0, 0, 0, 0, 0, 0, 0, 0, 0, 1, 0,\ldots.$$
Thus both sequences $a_n$ and $a_n*$ have the same Hankel transform. However, the Hankel transform of $a_n$ corresponds to the
pattern sequence $0,1,3,5,11,\ldots$ while that of $a_n*$ has distribution $0,1,1,3,5,11,\ldots$. To distinguish the Hankel transforms, we could accompany each ``continued-fraction-derived'' $\{-1,0,1\}$ Hankel transform with its Hankel pattern sequence (or equivalently its CF power sequence). Using the Hankel pattern sequence for this purpose is equivalent to assigning multiplicities to the non-zero elements of the Hankel transform. Thus in the case of the Hankel transform of $a_n*$, we assign a multiplicity of two to the second $1$ in the Hankel transform, corresponding to the repetition of element $J_1=J_2$. To emphasize this, we sometimes write this Hankel transform as
$$1, 1_2, 0, -1, 0, 1, 0, 0, 0, 0, 0, -1, 0, 0, 0, 0, 0, 0, 0, 0, 0, 1, 0,\ldots, $$ where un-indexed numbers have multiplicity $1$.

\end{example}

In subsequent examples we allow for repeated elements in the Hankel pattern sequence. The next case is an extreme case, where all elements (bar the first) are duplicated.

\begin{example} The Catalan numbers $C_n=\frac{1}{n+1} \binom{2n}{n}$ \seqnum{A000108} with generating function $c(x)$ given by
$$c(x)=\frac{1-\sqrt{1-4x}}{2x}= \cfrac{1}{1-
\cfrac{x}{1-
\cfrac{x}{1-
\cfrac{x}{1-\cdots}}}}
$$ is well known to have Hankel transform $1,1,1,\ldots$ \cite{Kratt}; that is, each element $h_n=1$. The CF power sequence for $C_n$ is also $1,1,1,\ldots$. This means that the Hankel pattern for $C_n$ is the sequence $$0,1,1,2,2,3,3,4,4,\ldots$$ with general term $\lfloor \frac{n+1}{2} \rfloor$. Thus the first term $h_0=1$ has multiplicity $1$ while all the other terms $h_n=1$ have multiplicity $2$.
Using the index notation, we could thus write the Hankel transform of $C_n$ as
$$1,1_2,1_2,1_2,1_2,1_2,1_2,1_2,1_2,1_2,\ldots.$$
\end{example}

\begin{example} We consider the CF power sequence $\{p_n\}$ given by $1,1,2,2,2,2,2,\ldots$. This corresponds to the Hankel pattern
$$0,1,2,3,4,5,6,\ldots.$$ The sequence $a_n$ with generating function
$$\cfrac{1}{1-
\cfrac{x}{1-
\cfrac{x}{1-
\cfrac{x^2}{1-
\cfrac{x^2}{1-
\cfrac{x^2}{1-
 \cdots}}}}}}$$
 has Hankel transform $1,1,1,1,\ldots$ where each ``$1$'' has multiplicity $1$. Thus we could have written this Hankel transform as
 $$1_1,1_1,1_1,1_1,1_1,1_1,1_1,1_1,\ldots.$$

\noindent The sequence $a_n$ begins
 $$1, 1, 2, 4, 9, 20, 46, 105, 243, 560, 1299, 3006, \ldots,$$
\noindent The generating function of this sequence is
$$\cfrac{1}{1-\cfrac{x}{1-x c(x^2)}}.$$

 \end{example}

\begin{example} It is not the case that every $\{-1,0,1\}$ sequence is the Hankel transform of an integer sequence. For example, the sequence $$1,0,1,0,1,0,1,0,1,\ldots $$ is not the Hankel transform of any integer (or real) sequence.
\begin{proof} (Somos)
Assume that the sequence $\alpha, \beta, \gamma, \delta, \epsilon,\ldots$ has Hankel transform $1,0,1,0,\ldots$.
We have $h_0=\alpha$ and hence $\alpha=1$. 

\noindent Then $h_1=\gamma-\beta^2=0$ and so $\gamma=\beta^2$, and so the sequence would start $1,\beta,\beta^2, \delta, \ldots$.

\noindent Then $h_2=\beta^6+2\beta^3 \delta-\delta^2$, and so we must have  
$$\beta^6+2\beta^3 \delta-\delta^2=1.$$ The solution of this quadratic in $\delta$ is 
$$ \delta=\beta^3-i \quad \text{or} \quad \delta=\beta^3+i.$$ 
\end{proof}

\end{example}

\begin{example} Consider the CF power sequence $p_n$ given by $1,1,3,3,3,3,3,\ldots$. This gives us the Hankel pattern sequence
$$0, 1, 3, 4, 6, 7, 9, 10, 12, 13, 15, \ldots.$$

The corresponding sequence begins
$$1, 1, 2, 4, 8, 17, 36, 76, 162, 345, 734, 1565, 3336, 7109, 15158, 32318, 68898, \ldots$$ with generating function  $$\cfrac{1}{1-\cfrac{x}{1-x c(x^3)}}.$$
\noindent Its Hankel transform is equal to the periodic sequence
$$1, 1, 0, -1, -1, 0, 1, 1, 0, -1, -1, 0, \ldots.$$
\noindent Here, the non-zero terms all have multiplicity one.
\end{example}
\section{A number theoretic example}
We let the Hankel  pattern sequence be the sequence $b_n=\binom{n+1}{2}$ of the triangular numbers. We find that the corresponding CF power sequence is the sequence $$1,1,3,5,7,9,11,\ldots$$ of extended odd numbers. The sequence generated by the continued fraction defined by the power sequence begins
$$ 1, 1, 2, 4, 8, 17, 36, 76, 161, 341, 723, 1533, 3250, 6891, 14611, 30980, 65688, 139281, \ldots $$ and has a Hankel $h_n$ transform
that begins
$$ 1, 1, 0, -1, 0, 0, 1, 0, 0, 0, 1, 0, 0, 0, 0, 1, 0, 0, 0, 0, 0, -1, 0, \ldots $$

In this case, the quantity
$$e_n=\sum_{k=0}^n (-1)^{n-k} h_n h_{n-k}$$ is of interest. This is the convolution of $h_n$ with $(-1)^n h_n$. The sequence
$e_{2n}$ begins $$1, -1, 2, 1, 0, 2, 1, 0, 0, 2, 1, 2, \ldots $$ and is related to the so-called ``eta quotients''.

If we now take the CF power sequence $p_n = \lfloor \frac{n+1}{2} \rfloor +0^n$, or $1,1,1,2,2,3,3,4,4,5,5,\ldots$, we obtain
the Hankel pattern sequence $$0, 1, 1, 3, 3, 6, 6, 10, 10, 15, 15, \ldots $$
The sequence $a_n$ now begins
$$1, 1, 2, 5, 13, 35, 95, 260, 713, 1959, 5386, 14815, 40759, 112151, 308609, 849240, 2337009, 6431246, \ldots $$ and has Hankel transform $$1, 1, 0, -1, 0, 0, 1, 0, 0, 0, 1, 0, 0, 0, 0, 1, 0, 0, 0,\ldots$$
\noindent Here, apart from the first term, the non-zero terms have multiplicity two.

\section{More variations on $1,3,5,\ldots$}
We start this section by noting that the sequence defined by the CF power sequence $1,3,5,7,\ldots$ is \seqnum{A143951}, which counts the  number of Dyck paths such that the area between the $x$-axis and the path is $n$ (Emeric Deutsch). The Hankel transform of this sequence begins 
$$1, 0, 0, -1, 0, 1, 0, 0, 0, 0, 1, 0, 0, 0, 1, 0, 0, 0, 0, 0, 0, -1,\ldots,$$ with Hankel pattern sequence 
$$b_n=(n^2+3n+1+(-1)^n)/2,$$ which is essentially \seqnum{A176222}. The sequence $b_n$ thus begins 
$$0, 3, 5, 10, 14, 21, 27, 36, 44, 55, 65,\ldots.$$

For the CF power sequence $1,1,3,3,5,5,7,7,\ldots$ we find that the sequence $b_n$, which begins 
$$0, 1, 3, 4, 8, 9, 15, 16, 24, 25, 35,\ldots,$$ satisfies
$$b_n=(2n^2+6n+1+(2n-1)(-1)^n)/8.$$ 
The corresponding sequence $a_n$ begins 
$$1, 1, 2, 4, 8, 17, 36, 76, 162, 345, 734, 1564, 3332, \ldots$$ with Hankel transform 
$$ 1, 1, 0, -1, -1, 0, 0, 0, -1, -1, 0, 0, 0, 0, 0, 1, 1, 0, 0, \ldots.$$ 

For the CF power sequence $1,1,1,3,3,3,5,5,5,7,7,7,\ldots$  we find that the sequence $b_n$, which begins 
$$0,1, 1, 3, 3, 5, 6, 8, 9, 12, 13, 16,\ldots,$$ has generating function 
$$\frac{x(1+x^2-x^3)}{(1-x)^3 (1+2x+2x^2+x^3)}.$$ The corresponding sequence $a_n$ begins 
$$ 1, 1, 2, 5, 13, 34, 90, 239, 635, 1689, 4494, 11958, 31823, 84692, 225396, \ldots$$ and has Hankel transform 
$$1, 1, 0, 0, -1, 0, 0, 1, 0, -1, 0, 0, 1, 0, -1, 0,\ldots. $$

\noindent Some slight variations on this last example are also of interest. For instance, the CF power sequence 
$$1, 2, 3, 3, 3, 5, 5, 5, 7, 7, 7, 9, 9, 9, 11, \ldots$$ corresponds to the Hankel pattern sequence $b_n$ which begins 
$$0, 2, 3, 5, 6, 10, 11, 15, 18, 22, 25, 31, 34, 40, 45,\ldots,$$ with generating function 
$$\frac{x(2+x-2x^3+x^4)}{(1-x)^3(1+2x+2x^2+x^3)}.$$ The Hankel transform of $a_n$ is then 
$$1, 0, -1, -1, 0, 1, 1, 0, 0, 0, 1, 1, 0, 0, 0, 1,\ldots.$$ 

\section{A pattern avoiding example} The CF power sequence $$1, 1, 1, 1, 2, 1, 1, 2, 1, 1, 2, 1, 1, 2, 1, 1, 2, 1, 1, 2, 1,\ldots $$ corresponds to the Hankel pattern sequence $$0, 1, 1, 2, 3, 3, 4, 5, 5, 6, 7, 7, 8, 9, 9, 10, 11, 11, 12, 13, 13, \ldots $$ which gives a $\{-1,0,1\}$ Hankel transform where odd indexed non-zero terms have multiplicity two, and the even-indexed non-zero terms have multiplicity one.
The corresponding sequence is \seqnum{A054391} \cite{Mansour_Shattuck}. It begins
$$1, 1, 2, 5, 14, 41, 123, 374, 1147, 3538, 10958, 34042, 105997, \ldots$$  and  has Hankel transform consisting of all $1$'s (with the multiplicities above). Thus we could write the Hankel transform as
$$1,1_2,1,1_2,1,1_2,1,\ldots.$$

We note that the generating function of this sequence is given by
$$\cfrac{1}{1-
\cfrac{x}{1-x u}},$$ where $u$ satisfies the equation
$$u=\cfrac{1}{1-\cfrac{x}{1-\cfrac{x}{1-x^2 u}}}.$$
Solving, we find that the generating function is equal to
$$\frac{1-3x-\sqrt{1-2x-3x^2}}{1-3x+2x^2+\sqrt{1-2x-3x^2}}, $$ confirming that this sequence is the same as \seqnum{A054391}.

\section{The Motzkin numbers} We consider the CF power sequence
$$ 1, 1, 2, 1, 1, 2, 1, 1, 2, 1, 1, 2, 1, 1, 2, 1, 1, 2, 1, 1, 2, \ldots.$$ This corresponds to the Hankel pattern sequence
$$0, 1, 2, 2, 3, 4, 4, 5, 6, 6, 7, 8, 8, 9, 10, 10, 11, 12, 12, 13, 14, \ldots $$ so once again we obtain a sequence with an
all $1$'s Hankel transform, but with the multiplicities indicated. Thus we get the Hankel transform
$$1,1,1_2,1,1_2,1,1_2,1,\ldots.$$

The generating function $u=g(x)$ satisfies the equation
$$u=\cfrac{1}{1-\cfrac{x}{1-\cfrac{x}{1-x^2 u}}},$$ which solves to give
$$g(x)=\frac{1-x-\sqrt{1-2x-3x^2}}{2x^2},$$ which coincides with the g.f. of the Motzkin numbers \seqnum{A001006}.
\section{Euler pentagonal numbers}
In this section, we take the Euler pentagonal numbers \seqnum{A001318} \cite{Somos} as the basis of the CF power sequence. Thus we define
$$p_n=\frac{6n^2+6n+1}{16}-\frac{(2n+1)(-1)^n}{16}+0^n$$ which begins
$$1, 1, 2, 5, 7, 12, 15, 22, 26, 35, 40,\ldots.$$ The corresponding Hankel pattern sequence is given by
$$b_n = \frac{1}{2}\lfloor \frac{n+1}{2} \rfloor \lfloor \frac{n+2}{2} \rfloor \lfloor \frac{n+3}{2} \rfloor,$$ which is
\seqnum{A028724}. This begins
$$0, 1, 2, 6, 9, 18, 24, 40, 50, 75, 90,\ldots .$$

The resulting sequence begins
$$1, 1, 2, 4, 9, 20, 45, 101, 227, 511, 1150, 2589, 5828, 13120, 29536, 66492, 149690,\ldots,$$ and has Hankel transform

$$1, 1, 1, 0, 0, 0, 1, 0, 0, -1, 0, 0, 0, 0, 0, 0, 0, 0, -1, 0, 0, \ldots$$

Note that if in the continued fraction we take the pattern of signs $-,-,+,+,-,-,+,+,\ldots$ then the resulting sequence begins
$$1, 1, 2, 4, 7, 12, 21, 37, 65, 115, 204, 361, 638, 1128, 1994, 3524, 6230, \ldots$$ and has a Hankel transform that begins
$$1, 1, -1, 0, 0, 0, -1, 0, 0, 1, 0, 0, 0, 0, 0, 0, 0, 0, 1, 0, 0,\ldots $$

\section{A Fibonacci distribution}
It is of interest to find a sequence whose Hankel transform has a $\{-1,0,1\}$ distribution that follows the Fibonacci numbers \seqnum{A000045}. Corresponding to the CF power sequence
$$1,1,1,1,2,3,5,8,13,21,34,55,89,\ldots $$ we get the Hankel pattern sequence
$$0,1,1,2,3,5,8,13,21,34,55,89,\ldots.$$
The sequence defined by the corresponding continued fraction begins
$$1, 1, 2, 5, 14, 41, 123, 373, 1137, 3475, 10634, 32562, 99738, 305546, 936108, 2868084,\ldots$$
and it has a Hankel transform that starts
$$1, 1, 1, 1, 0, -1, 0, 0, 1, 0, 0, 0, 0, 1, 0, 0, 0, 0, \ldots.$$

\section{A gap Hankel transform}
Consider the CF power sequence
$$1,3,4,2,2,2,2,2,\ldots.$$
The corresponding Hankel pattern sequence is
$$0, 3, 4, 5, 6, 7, 8, 9, 10, 11, 12,\ldots.$$
Thus there will be a 'gap' at positions $1$ and $2$ (filled by zeros). In fact, the sequence $a_n$ corresponding to this CF power sequence begins
$$1, 1, 1, 1, 2, 3, 4, 6, 10, 15, 23, 36, 58, 90, 145, 230, 377, 601, 1000,\ldots$$
and has Hankel transform $h_n$ given by
$$1,0,0,-1,-1,-1,-1,-1,-1,-1,-1,-1,-1,-1,-1,\ldots.$$
\noindent The sequence $a_n$ has some interesting Hankel properties. The sequence $a_{n+1}$ has Hankel transform given by
$$1,0,-1,0,1,0,-1,0,1,0,\ldots.$$ \noindent The twice shifted sequence $a_{n+2}$ has Hankel transform
$$1, 1, -1, -1, -2, -2, -3, -3, -4, -4, \ldots,$$ while the three-times shifted sequence $a_{n+3}$ has Hankel transform
$$1, -1, -1, 1, 4, -1, -9, 1, 16, \ldots.$$ If we prepend $1,1$ to the sequence then the new sequence has a Hankel transform
of
$$1,0,0,0,0,1,1,2,2,3,3,\ldots.$$
Other sequences with ``gap'' Hankel transforms can be built using the template
$$1,r,r+1,2,2,2,2,2,2,2,\ldots$$ for the CF power sequence, since this maps to the pattern sequence
$$0, r, r + 1, r + 2, r + 3, r + 4, r + 5, r + 6, r + 7, r + 8, r + 9,\ldots.$$

\section{Conclusion} The foregoing shows that $\{-1,0,1\}$ Hankel transforms are objects worthy of study. The notion of multiplicity seems important, particularly as a way of further characterizing certain sequences that would otherwise have the ``same'' Hankel transforms. Thus the Catalan numbers $C_n$ have Hankel transform $1,1,1,1,\ldots$ but with pattern sequence $0,1,1,2,2,3,3,\ldots,$ while the Motzkin numbers $M_n$ also have Hankel transform $1,1,1,1,\ldots$ but with pattern sequence $0,1,2,2,3,4,4,\ldots$.

\section{Acknowledgement} The material in this note arose out of an exchange of emails with Michael Somos, who originally posed the question of the origins of $\{-1,0,1\}$ Hankel transforms. I am extremely grateful to him for the time spent on this exchange.

\bigskip \hrule \bigskip 

\noindent 2010 {\it Mathematics Subject Classification}: Primary 11B83; Secondary 05A15, 11J70. 
\\
\noindent \emph{Keywords:}
Permutation, integer sequence, Hankel determinant, Hankel
transform, continued fraction

 \bigskip \hrule \bigskip \noindent Concerns sequences

\seqnum{A000045},
\seqnum{A000108},
\seqnum{A001006},
\seqnum{A001045},
\seqnum{A001318},
\seqnum{A028724},
\seqnum{A054391},
\seqnum{A143951}.


\begin{thebibliography}{13}

\bibitem{Kratt}
C. Krattenthaler, Advanced determinant calculus: A
complement, \emph{Lin. Alg. Appl.} {\bf 411} (2005)
68–-166.

\bibitem{Layman} J. W. Layman,
 The Hankel transform and some
    of
    its properties, \emph{J. Integer Seq.} {\bf
    4} (2001), \href{http://www.cs.uwaterloo.ca/journals/JIS/VOL4/LAYMAN/hankel.html}{Article 01.1.5}.

\bibitem{Mansour_Shattuck}
T. Mansour, M. Shattuck, Pattern avoiding partitions, sequence \seqnum{A054391} and the kernel method, \emph{Applications and Applied Mathematics} \textbf{6} (2011), 397--411.

\bibitem{SL1} N. J. A.~Sloane, \emph{The
On-Line Encyclopedia of Integer Sequences}. Published electronically
at \href{http://oeis.org}{http://oeis.org}, 2012.

\bibitem{Somos}
M. Somos, A Multisection of $q$-Series, http://cis.csuohio.edu/~somos/multiq.pdf

\end{thebibliography}
\end{document}